\documentclass[a4paper,10pt,conference]{IEEEconf}
\usepackage{cite}
\usepackage{amsmath}
\usepackage{amssymb}
\usepackage{graphicx}
\usepackage{psfrag}
\begin{document}
\newtheorem{theorem}{Theorem}[section]
\newtheorem{definition}{Definition}[section]

\title {On the Physical Realizability of a Class of Nonlinear Quantum Systems}
\author{Aline I. Maalouf and Ian R. Petersen \\
School of Engineering and Information Technology,\\ University of New South Wales
at the Australian Defence Force
Academy,\\ Canberra, ACT 2600 \tt\small a.maalouf@adfa.edu.au \tt\small i.r.petersen@gmail.com
\thanks{This work was completed with the support of a University
of New South Wales Postgraduate Award and the Australian Research Council.}}
\maketitle
\thispagestyle{empty}
\pagestyle{empty}
\maketitle
\begin{abstract}
In this paper, the physical realizability property is investigated for a class of nonlinear quantum systems. This property determines whether a given set of nonlinear quantum stochastic differential equations corresponds to a physical nonlinear quantum system satisfying the laws of quantum mechanics.
\end{abstract}
\section{Introduction}
In the physics literature, methods have been developed to model a wide range of open quantum systems, such as those encountered in quantum optics, within the framework of quantum stochastic differential equations \cite{Gardiner85}, \cite{Gardiner00}, \cite{Hudson84}. In this context, quantum noise is used to represent the influence of large heat baths and boson fields, including optical and phonon fields, from which completely positive maps, Lindblad generators and master equations are obtained by taking expectations \cite{James10}.
We may then distinguish two types of quantum stochastic differential equations, namely, linear and nonlinear. Linear quantum systems can be described by linear quantum stochastic differential equations and arise mostly in the area of quantum optics; e.g., see \cite{Walls08}, \cite{Gardiner00} and \cite{Bachor04}.
An important class of linear quantum stochastic models describe the Heisenberg evolution of the annihilation and creation operators of several independent open quantum harmonic oscillators that are coupled to external coherent bosonic fields, such as coherent laser beams; e.g., see \cite{Wiseman10}, \cite{Walls08} and \cite{Gardiner00}. A special class of these linear quantum systems is driven by quantum Wiener processes as in \cite{James08} where physical realizability conditions are developed to determine when the linear quantum system under consideration can be regarded as a representation of a linear quantum harmonic oscillator. This imposes some restrictions on the matrices describing the linear quantum model. This notion of physical realizability has been further investigated in \cite{Shaiju09} where the authors prove the equivalence between the algebraic conditions for physical realizability obtained in \cite{James08} and a frequency domain condition that an associated linear system is $(J, J)$-unitary. In fact, \cite{Shaiju09} extends the frequency domain physical realizability results of \cite{Maalouf09}, \cite{Maalouf11-1} and \cite{Maalouf11-2} which apply to linear quantum systems described purely in terms of annihilation operators only to a more general class of linear quantum systems, which are described in terms of both the annihilation and creation operators. More explicitly, in \cite{Maalouf09}, \cite{Maalouf11-1} and \cite{Maalouf11-2}, the physical realizability conditions developed for annihilation-operator linear quantum systems have been related to the lossless property of linear systems and \cite{Shaiju09} generalizes this result to relate the corresponding physical realizability conditions of annihilation and creation operators linear quantum systems to the property of $(J, J)$-unitary.
However, the question of addressing the physical realizability conditions for nonlinear quantum systems described by nonlinear quantum stochastic differential equations is still open. And in this paper, we restrict attention to a class of nonlinear quantum systems which is a generalization of the annihilation operator only linear quantum systems considered in \cite{Maalouf09}, \cite{Maalouf11-1} and \cite{Maalouf11-2}.

The paper is organized as follows: Section II defines the class of nonlinear quantum systems under consideration. Section III develops conditions under which the commutation relations for the nonlinear quantum systems are preserved. Section IV defines physical realizability in terms of nonlinear open quantum harmonic oscillators. Section V provides an example to illustrate the theory developed. Section VI concludes the paper.

\section{The Class of Nonlinear Quantum Systems}

We consider an open quantum system G with physical variable space $\mathcal{A}_G$ consisting of operators defined on an underlying Hilbert space $\mathcal{H}_G$. The self-energy of this system is described by a Hamiltonian $\mathcal{H}\in \mathcal{A}_G$. The system is driven by $m$ field channels given by the quantum stochastic process $W$

\begin{equation}\nonumber
W= \left( {\begin{array}{*{20}{c}}
{\begin{array}{*{20}{c}}
{{W_1}}\\
 \vdots \\
{{W_m}}
\end{array}}
\end{array}} \right).
\end{equation}
These describe the annihilation of photons in the field channels and are operators on a Hilbert space $F$, with associated variable space $\mathcal{F}$. In that case, $F$ is the Hilbert space defining an indefinite number of quanta (called a Fock space \cite{Parthasarathy92}), and $\mathcal{F}$ is the space of operators over this space.

We assume the process $W$ is canonical, meaning that we have the following second order Ito products:
\begin{eqnarray}\nonumber
dW_k(t)dW_l(t)^*&=&\delta_{kl}dt;\\\nonumber
dW_k(t)^*dW_l(t)&=&0;\\\nonumber
dW_k(t)dW_l(t)&=&0;\\\nonumber
dW_k(t)^*dW_l(t)^*&=&0
\end{eqnarray}
where $W_k(t)^*$ is the operator adjoint of $W_k(t)$ which is defined on the same Fock space.

The system is coupled to the field through a scattering matrix $S=I$ and a coupling vector of operators $L$ given by

\begin{equation}\nonumber
L=\left( {\begin{array}{*{20}{c}}
{\begin{array}{*{20}{c}}
{{L_1}}\\
 \vdots \\
{{L_m}}
\end{array}}
\end{array}} \right)
\end{equation}
where $L_j \in \mathcal{A}_G$.

The notation $G=(S,L,H)$ is used to indicate an open quantum system specified by the parameters $S, L$ and $H$ where $H$ represents the Hamiltonian of the system. For $S=I$, the Schrodinger equation is given by
\begin{equation}\label{schrodinger}
dU(t)= \left\{dW^\dagger L - L^\dagger dW -\frac{1}{2}L^\dagger L dt-iHdt\right\}U(t)
\end{equation}
with the initial condition $U(0)=I$. Equation (\ref{schrodinger}) determines the unitary motion of the system in accordance of the fundamental laws of quantum mechanics. Note that the notation $^\dagger$ refers to the Hilbert space adjoint.

Given a system annihilation operator $a_l \in \mathcal{A}_G$, its Heisenberg evolution is defined by $a_l(t)=j_t(a_l)=U(t)^* a_l U(t)$ and satisfies
\begin{eqnarray}\label{annihi_eq}
da_l(t)&=&\left(\mathcal{L}_L(t)(a_l(t))-i\left[a_l(t),H(t)\right]\right)\\\nonumber
&& + \left[L(t)^\dagger,a_l(t)\right]dW(t)
\end{eqnarray}
where the notation $\left[A,B\right]=AB-BA$ denotes the commutator of two operators $A$ and $B$.

In equation (\ref{annihi_eq}), all operators evolve unitarily and the notation $\mathcal{L}_L(a_l)$ refers to
\begin{equation}\label{gen_l}
\mathcal{L}_L(a_l)= \frac{1}{2}L^\dagger \left[a_l,L\right]+ \frac{1}{2} \left[L^\dagger,a_l\right]L.
\end{equation}
In this paper, we restrict attention to annihilation only coupling operators; i.e, the coupling operator $L(t)$ depends only on the annihilation operator $a_l(t)$ and not the creation operator $a_l(t)^*$. Therefore, $L(t)$ satisfies the following commutation property: $\left[a_l,L\right]=0$ and then
\begin{equation}\label{gen_l_a}
\mathcal{L}_L(a_l)= \frac{1}{2} \left[L^\dagger,a_l\right]L.
\end{equation}
The generator of the system $G$ is given by
\begin{equation}\label{gen_g}
\mathcal{G}_G(a_l)=-i\left[a_l,H\right]+\mathcal{L}_L(a_l).
\end{equation}
For the case of having $n$ annihilation operators $a_1,\ldots,a_n$, we define
\begin{equation}\nonumber
a= \left( {\begin{array}{*{20}{c}}
{\begin{array}{*{20}{c}}
{{a_1}}\\
 \vdots \\
{{a_n}}
\end{array}}
\end{array}} \right),
\end{equation}
\begin{equation}\label{genl_i}
\mathcal{L}_{{L}_i}(a)=  \frac{1}{2} \left[L^\dagger,a_i\right]L
\end{equation}
and
\begin{equation}\label{geng_i}
\mathcal{G}_{G_i}(a)=-i\left[a_i,H\right]+\mathcal{L}_{L_i}(a).
\end{equation}
Therefore, we can write
\begin{eqnarray}\nonumber \label{annihieq_i}
da_i(t)&=&\left(\mathcal{L}_{L_i}(t)(a_i(t))-i\left[a_i(t),H(t)\right]\right)\\\nonumber
&& + \left[L(t)^\dagger,a_i(t)\right]dW(t)
\end{eqnarray}
and
\begin{eqnarray}\nonumber \label{annihieq_a}
da(t)&=&\left(\mathcal{L}_{L}(t)(a(t))-i\left[a(t),H(t)\right]\right)\\\nonumber
&& + \left[L(t)^\dagger,a(t)\right]dW(t)
\end{eqnarray}
where
\begin{equation}\label{L_a}
\mathcal{L}_{L}(a)= \left( {\begin{array}{*{20}{c}}
{\begin{array}{*{20}{c}}
{{\mathcal{L}_{{L}_1}(a)}}\\
 \vdots \\
{{\mathcal{L}_{{L}_n}(a)}}
\end{array}}
\end{array}} \right)
\end{equation}
and
\begin{equation}\label{G_a}
\left[a,H\right]= \left( {\begin{array}{*{20}{c}}
{\begin{array}{*{20}{c}}
{{\left[a_1,H\right]}}\\
 \vdots \\
{{\left[a_n,H\right]}}
\end{array}}
\end{array}} \right).
\end{equation}
Let
\begin{eqnarray}\nonumber\label{At}
A(a(t), a(t)^\dagger)&=&\left({\begin{array}{*{20}{c}}
{\begin{array}{*{20}{c}}
{{A_1(a(t),a(t)^\dagger )}}\\
 \vdots \\
{{A_n(a(t), a(t)^\dagger)}}
\end{array}}
\end{array}}\right)\\\nonumber&=& \left( {\begin{array}{*{20}{c}}
{\begin{array}{*{20}{c}}
{{\mathcal{L}_{{L}_1}(a)}}\\
 \vdots \\
{{\mathcal{L}_{{L}_n}(a)}}
\end{array}}
\end{array}} \right)-i\left( {\begin{array}{*{20}{c}}
{\begin{array}{*{20}{c}}
{{\left[a_1,H\right]}}\\
 \vdots \\
{{\left[a_n,H\right]}}
\end{array}}
\end{array}} \right)\\\nonumber &=& \mathcal{L}_{L}(a)-i\left[a,H\right]\\ &=&\mathcal{G}_G(a)
\end{eqnarray}
and
\begin{eqnarray}\label{Bt}
B(a(t),a(t)^\dagger)&=&\left({\begin{array}{*{20}{c}}
{\begin{array}{*{20}{c}}
{{B_1(a(t),a(t)^\dagger)}}\\
 \vdots \\
{{B_n(a(t),a(t)^\dagger)}}
\end{array}}
\end{array}}\right)\\\nonumber &=& \left({\begin{array}{*{20}{c}}
{\begin{array}{*{20}{c}}
{{\left[L^\dagger, a_1\right]}}\\
 \vdots \\
{{\left[L^\dagger, a_n\right]}}
\end{array}}
\end{array}}\right)\\&=& \left[L^\dagger, a\right].
\end{eqnarray}
Hence,
\begin{eqnarray}\nonumber
da(t)&=&A(a(t),a(t)^\dagger)dt+B(a(t),a(t)^\dagger)dW(t)\\\nonumber
da(t)^*&=&A(a(t),a(t)^\dagger)^*dt+B(a(t),a(t)^\dagger)^*dW(t)^*.
\end{eqnarray}
Note that for the case of matrices, the notations $^*$ and $^\dagger$ refer respectively to the complex conjugate and the complex conjugate transpose of the matrix in question.
\begin{equation}\label{at}\end{equation}
On the other hand, the components of the output fields are defined by $y(t)=j_t(W(t))=U(t)^*W(t)U(t)$ and satisfy the nonlinear quantum stochastic differential equations
\begin{eqnarray}\nonumber
dy(t)&=&C(a(t))dt+D(t)dW(t)\\\nonumber
dy(t)^*&=&C(a(t))^*dt+D(t)^*dW(t)^*
\end{eqnarray}
\begin{equation}\label{yt}\end{equation}
where $C(a(t))=L(t)$, $C(a(t))^*=L(t)^*$, $D(t)=I$ and $D(t)^*=I$.
Hence, the system $G$ can be described by the following nonlinear quantum stochastic differential equations
\begin{eqnarray}\nonumber\label{sys_a}
d\bar a(t)&=&\bar A(a(t),a(t)^\dagger)dt+\bar B(a(t),a(t)^\dagger)d\bar W(t);\\
d\bar y(t)&=&\bar C(a(t),a(t)^\dagger)dt+\bar D(a(t),a(t)^\dagger)d\bar W(t)
\end{eqnarray}
where $\bar a (t)=\left[ {\begin{array}{*{20}{c}}
{a}\\
{a^*}
\end{array}} \right]$, $\bar A(a,{a^\dag }) = \left[ {\begin{array}{*{20}{c}}
{A(a,{a^\dag })}\\
{A{{(a,{a^\dag })}^*}}
\end{array}} \right]$, $\bar C(a,{a^\dag }) = \left[ {\begin{array}{*{20}{c}}
{C(a)}\\
{C{{(a)}^*}}
\end{array}} \right]$, $\bar B(a,{a^\dag }) = \left[ {\begin{array}{*{20}{c}}
{B(a,a^\dagger)}&0\\
0&{B{{(a,a^\dagger)}^*}}
\end{array}} \right]$, $\bar D(t) = \left[ {\begin{array}{*{20}{c}}
{D(t)}&0\\
0&{D{{(t)}^*}}
\end{array}} \right]$ and $d\bar W(t)=\left[ {\begin{array}{*{20}{c}}
{dW(t)}\\
{dW(t)^*}
\end{array}} \right]$.
\begin{definition}\label{def4}
The class of nonlinear quantum system we consider in this paper are nonlinear quantum systems that can be represented by the QSDES (\ref{sys_a}) such that the matrices $A(a,a^\dagger)$ and $ C(a)$ satisfy $\left[C(a), a^T\right]=0$, $\left[A(a,a^\dagger), a^T\right]=-\left[a, A(a,a^\dagger)^T\right]$ and
\begin{eqnarray}\nonumber \label{prop_matrices}
A_i(a,a^\dagger)&=& \sum\limits_{{k_i} = 0}^{{m_{{k_i}}}} {\sum\limits_{{h_i} = 0}^{{m_{{h_i}}}} {\sum\limits_{l = 1}^n {\sum\limits_{p = 1}^n {{\alpha _{{plk_ih_i}} } }}}} a_p^{{k_i}}{(a_l^*)^{{h_i}}}\\
C_v(a)&=& \sum\limits_{{k_v} = 0}^{{m_{{k_v}}}}\sum\limits_{p = 1}^n \beta_{k_vp}a_p^{k_v}
\end{eqnarray}
for $i=1,\ldots,n$ and $v=1,\ldots, m$ with $m_{k_i}$, $m_{h_i}$, $m_{k_v}$, $\beta_{k_vp}$ and $\alpha _{{plk_i}{h_i}}$ integers.

Moreover, the matrices $A(a,a^\dagger)$, $B(a,a^\dagger)$ and $C(a)$ satisfy the following property:
\begin{eqnarray}\nonumber \label{prop_matrices_2}
&&\frac{1}{{{2{\bar n}}}}\left[ {\bar A{{(a,{a^\dag })}^\dag }{{\bar \theta }^{ - 1}}\bar a,\bar a} \right] - \frac{1}{{{2{\bar n}}}}\left[ {{{\bar a}^\dag }{{\bar \theta }^{ - 1}}\bar A(a,{a^\dag }),\bar a} \right]=\\&& \qquad \bar A(a,{a^\dag }) - \frac{1}{2}\bar B(a,{a^\dag })\bar C(a,a^\dagger)
\end{eqnarray}
where $\bar n= {\left. {\mathop {\sup }\limits_{1 \le i \le n} \left( {{k_i} + {h_i}} \right)} \right|_{{\alpha _{{plk_ih_i}}}\ne 0}}$, $\bar \Theta= \left[ {\begin{array}{*{20}{c}}
\Theta &0\\
0&{{\Theta ^*}}
\end{array}} \right]$ and $\Theta = \left[a,a^\dagger\right]$. For more information, on $\Theta$, refer to Section III below.
\end{definition}

In the sequel, we assume that $dW(t)$ admits the following
decomposition:
\begin{equation}\nonumber
 dW (t) = \beta _{w} (t)dt + d\tilde w(t)
\end{equation}
where $\tilde w(t)$ is the noise part of $W(t)$ and $\beta_{w}(t)$
is an adapted process (see \cite{Parthasarathy92} and \cite{Hudson84}).
The noise $\tilde w(t)$ is a vector of quantum Weiner processes with Ito table
\begin{equation}\label{wIto}
d\tilde w(t)d\tilde w(t)^\dagger = F_{\tilde{w}} dt
\end{equation}
(see \cite{Parthasarathy92}) where $F_{\tilde{w}}$ is a non-negative definite Hermitian matrix. Here, the notation $^\dagger$ represents the adjoint transpose of a vector of operators. Also, we assume the following commutation relations hold for the noise components:
\begin{equation}
\label{noisecom}
\left[{d\tilde w(t),d\tilde w(t)^{\dagger}} \right] \triangleq d\tilde
w (t)d\tilde w(t)^{\dagger} - (d\tilde w(t)^* d\tilde
w(t)^T)^T = T_w dt.
\end{equation}
Here $T_w$ is a Hermitian matrix and the notation $^T$ denotes the transpose of a vector or matrix of operators. The noise processes can be represented as operators on an appropriate Fock
space; for more details, see \cite{Parthasarathy92}.

The process $\beta_{w}(t)$ represents variables of other systems
which may be passed to the system (\ref{sys_a}) via an interaction. Therefore,
it is required that $\beta_{w}(0)$ be an operator on a Hilbert
space distinct from that of $a_0$ and the noise processes. We also
assume that $\beta_{w}(t)$ commutes with $a(t)$ for all $t\geq
0$. Moreover, since $\beta_{w}(t)$ is an adapted process, we note that $\beta_{w}(t)$ also
commutes with $d\tilde{w}(t)$ for all $t\geq 0$. We have also that $\left[a, dW^\dagger \right]=0$.

\section{Commutation Relations}
For the nonlinear quantum system (\ref{sys_a}), the initial system variables $a(0)=a_{0}$ consist of operators satisfying the commutation relations
\begin{eqnarray}\nonumber
\label{relations}
\left[ {a_j (0),a_k^{*}(0)} \right] &=& \Theta _{jk} ,\qquad j,k = 1,
\ldots ,n.\\
\left[ {\bar a_j (0),\bar a_k^{*}(0)} \right] &=& \bar \Theta _{jk} ,\qquad j,k = 1,
\ldots ,2n.
\end{eqnarray}
Here, the commutator is defined by $\left[{a_j,a_k^{*}} \right]
\triangleq a_ja_k^{*}  - a_{k}^* a_j = \Theta_{jk}$ where $\Theta$ is a complex matrix with elements $\Theta_{jk}$.
With $a^T = (a_1,\ldots, a_n)$, the  relations (\ref{relations}) can  be written as
\begin{eqnarray}\nonumber \label{thetaandbar}
\left[ {a,a^\dagger } \right] & \triangleq & aa^\dagger  - (a^*a^T
)^T= \Theta\\
\left[\bar a, \bar a^\dagger\right]& \triangleq & \left[ {\begin{array}{*{20}{c}}
\Theta &0\\
0&{{\Theta ^*}}
\end{array}} \right]=\bar \Theta.
\end{eqnarray}
\subsection{Preservation of the Commutation Relations}
The following theorem provides an algebraic characterization of
when the nonlinear quantum system (\ref{sys_a}) preserves the commutation
relations (\ref{relations}) as time evolves.
\begin{theorem}\label{thm1}
For the nonlinear quantum system (\ref{sys_a}),
we have that $[\bar a_{p}(0), \bar a_{q}^{*}(0)]=\bar \Theta_{pq}$ implies
$[\bar a_p(t),\bar a_q^{*}(t)]=\bar \Theta_{pq}$ for all $t\geq 0$ with
$p,q=1\ldots, 2n$ if and only if
\begin{eqnarray}\nonumber \label{preservation}
&& \left[ \bar A(\bar a, \bar a^\dagger),\bar a^\dagger \right]+\left[\bar a, \bar A(\bar a,\bar a^\dagger)^\dagger \right]\\\nonumber && \quad + \bar B(\bar a, \bar a^\dagger)T_{ \bar w}\bar B(\bar a,\bar a^\dagger)^\dagger =0;\\\nonumber
&& \left[ \bar B(\bar a,\bar a^\dagger),\bar a^\dagger\right] = 0 \qquad \mbox{and}\\
&& \left[ \bar a, \bar B(\bar a,\bar a^\dagger)^\dagger \right]=0
\end{eqnarray}
with $T_{\bar w} = \left[ {\begin{array}{*{20}{c}}
{T_w}&0\\
0&{T_w^*}
\end{array}} \right]$.
\end{theorem}
\section{Physical Realizability And The Nonlinear Quantum Harmonic Oscillator}
In this section, we consider conditions under which a nonlinear quantum system of the form (\ref{sys_a}) corresponds to a nonlinear open quantum harmonic oscillator. Such a system is said to be physically realizable.
The class of nonlinear open quantum harmonic oscillators under consideration are defined by the nonlinear coupling operator $ \bar L $ and a nonlinear Hamiltonian $ \bar H$.  To derive a nonlinear system of the form (\ref{sys_a}) from a nonlinear open quantum harmonic oscillator defined by $\bar L$ and $\bar H$, we proceed as in Section II to get the following definition.
\begin{definition}\label{def1}
A nonlinear quantum system of the form (\ref{sys_a}) is said to be physically realizable if it is a representation of a nonlinear open quantum harmonic oscillator defined by a nonlinear coupling operator $\bar L$ and a nonlinear Hamiltonian $\bar H$, i.e., if there exist $\bar H$ and $\bar L$ such that the matrices $\bar A(a,a^\dagger), \bar B(a,a^\dagger), \bar C(a,a^\dagger)$ and $\bar D(a,a^\dagger)$ satisfy the following equations
\begin{eqnarray}\label{hamphys}\nonumber
\bar A(a,a^\dagger)&=&\frac{1}{2}\left[\bar L^\dagger, \bar a\right]\bar J \bar L + i\left[\bar H, \bar a\right] ;\\\nonumber
\bar B(a,a^\dagger)&=&\left[\bar L^\dagger, \bar a\right]\bar J;\\\nonumber
\bar C(a,a^\dagger)&=&\bar L;\\
\bar D(a,a^\dagger)&=&I;
\end{eqnarray}
where $\bar J =\left[ {\begin{array}{*{20}{c}}
I&{0}\\
{0}&-I
\end{array}} \right].$
\end{definition}
The following theorem provides necessary and sufficient conditions for a nonlinear quantum system of the form (\ref{sys_a}) to be physically realizable.
\begin{theorem}\label{thm2}
A nonlinear quantum system of the form (\ref{sys_a}) is  physically realizable with $T_{\bar w}= \bar J=\left[ {\begin{array}{*{20}{c}}
{I}&0\\
0&{-I}
\end{array}} \right] $ if and only if
\begin{eqnarray}\nonumber \label{phys_cond1}
&&\left[ \bar A( a,  a^\dagger),\bar a^\dagger \right]+\left[\bar a, \bar A( a, a^\dagger)^\dagger \right]\\\nonumber && \quad + \bar B( a,  a^\dagger)\bar J \bar B( a, a^\dagger)^\dagger =0;\\\nonumber
&& \left[ \bar B( a, a^\dagger),\bar a^\dagger \right] = 0;\\\nonumber
&& \left[ \bar a, \bar B( a, a^\dagger)^\dagger \right]=0;\\\nonumber
&& \bar B( a, a^\dagger)= \left[\bar C(a,a^\dagger)^\dagger,\bar a \right]\bar J\qquad \mbox{and}\\
&&\bar D( a,  a^\dagger)=I.
\end{eqnarray}
In this case, the corresponding nonlinear self-adjoint Hamiltonian $\bar H$ is given by
\begin{equation}\label{Ham}
\bar H=\frac{i}{2\bar n}\left({\bar a^\dag }{\bar \theta ^{ - 1}}\bar A(a,{a^\dag }) - \bar A{(a,{a^\dag })^\dag }{\bar \theta ^{ - 1}}\bar a\right)
\end{equation}
where $\bar n$ is as defined in Definition \ref{def4} and the corresponding coupling operator is $\bar L=\bar C(a,a^\dagger)$.
\end{theorem}
\section{Example}
We consider the application of the physical realizability conditions developed in this paper to an optical squeezer. A pulsed laser is used to drive the optical squeezer. The use of such a pulsed laser provides many advantages in an experimental situation; e.g, see \cite{Bachor04}.

The squeezer under consideration is implemented as an OPO driven by two time-varying optical fields; see Figure 1. A full description of the apparatus of Figure 1 is described in \cite{Hassen10}.

\begin{figure}
\includegraphics[width=8.5cm]{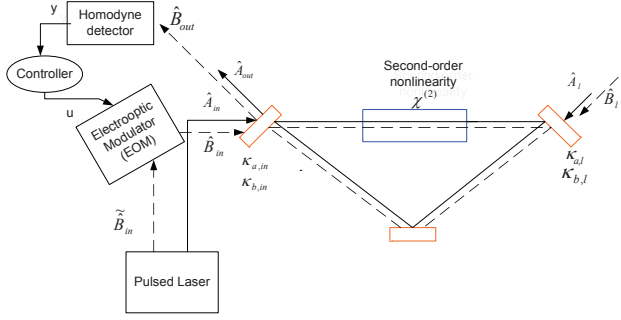}
\caption{Schematic for the control of an optical parametric oscillator.}
\end{figure}

\subsection{Modeling}

The dynamics of the OPO system can be described as in \cite{Hassen10} but with time-varying matrices and fields reflecting the fact that a pulsed laser is used to drive the experiment.
The dynamics of the optical subsystem contain nonlinear terms.
The quantities $A_{in}$ and $B_{in}$ are the input fields which set-up the
fundamental and second-harmonic intracavity fields $a_1$ and $a_2$. $\kappa_{a_{{in}}}$ and $\kappa_{b_{in}}$ are the loss rates of the input/output mirrors for the $a_1$ and $a_2$ fields respectively. The
parameters $\kappa_{a_{l}}$ and $\kappa_{b_l}$  are the internal loss rates for the corresponding two fields.
Also, $\kappa_{a_1}= \kappa_{a_{{in}}}+ \kappa_{a_{l}}$ and $\kappa_{a_2}= \kappa_{b_{{in}}}+ \kappa_{b_{l}}$ are the associated total resonator decay
rates. Furthermore, $\delta \hat A_l$ and $\delta \hat B_l$ represent the vacuum fields due to the internal
losses. The output fields are given by $\hat A_{out}$ and $\hat B_{out}$.

The Heisenberg equations of motion are given by (see \cite{Hassen10}):
\begin{eqnarray}\nonumber
d{a_1}(t)& =&A_1(a,a^\dagger)(t) dt + B_1(t)dw_1(t) \\\nonumber
d{a_2}(t) &=& A_2(a,a^\dagger)(t)dt + B_2(t) dw_2(t).
\end{eqnarray}
In that case,
$A_1(a,a^\dagger)= - {\kappa _{{a_1}}}{a_1} - 2{\chi ^{(2)}}a_1^*{a_2}$, $B_1=-\sqrt {2{\kappa _{{a_1}}}}$, $A_2(a,a^\dagger)=- {\kappa _{{a_2}}}{a_2} +{\chi ^{(2)}}a_1^2$, $B_2=-\sqrt {2{\kappa _{{a_2}}}}$, $dw_1(t)=-d{\hat A_{in}}$ and $dw_2(t)=-d{\hat B_{in}}$.
$a_1$ and $a_2$ satisfy the following commutation relations:
\begin{eqnarray}\nonumber
\left[ {{a_1},a_1^*} \right] &= &1;\\\nonumber
\left[ {{a_2},a_2^*} \right] &=&1 ;\\\nonumber
\left[ {{a_1},a_2^*} \right] &=&0; \\\nonumber
\left[ {{a_2},a_1^*} \right] &=&0;\\\nonumber
\left[ {{a_1},a_2} \right] &=&0.
\end{eqnarray}
Let $a = \left[ {\begin{array}{*{20}{c}}
{{a_1}}\\
{{a_2}}
\end{array}} \right]$, $A(a,{a^\dag }) = \left[ {\begin{array}{*{20}{c}}
{{A_1}(a,{a^\dag })}\\
{{A_2}(a,{a^\dag })}
\end{array}} \right]$, $B = \left[ {\begin{array}{*{20}{c}}
{{B_1}}&0\\
0&{{B_2}}
\end{array}} \right]$, $dw= \left[ {\begin{array}{*{20}{c}}
{d{w_1}}\\
{d{w_2}}
\end{array}} \right]$, $\bar a= \left[ {\begin{array}{*{20}{c}}
a\\
{  {a^*}}
\end{array}} \right]$, $\bar A(a,a^\dagger)= \left[ {\begin{array}{*{20}{c}}
{A(a,{a^\dag })}\\
{  A{{(a,{a^\dag })}^*}}
\end{array}} \right]$, $\bar B\left[ {\begin{array}{*{20}{c}}
{B}&0\\
0&{ B}^*
\end{array}} \right]$ and $d\bar w = \left[ {\begin{array}{*{20}{c}}
{dw}\\
{d{w^*}}
\end{array}} \right]$.

On the other hand, we define
\begin{equation}\nonumber
dy(t)=C(a)(t)dt+D(t)dw(t)
\end{equation}
where $C(a)= \left[ {\begin{array}{*{20}{c}}
{\sqrt {2{\kappa _{{a_1}}}} {a_1}}\\
{\sqrt {2{\kappa _{{a_2}}}} {a_2}}
\end{array}} \right]$ and $D(t)=\left[ {\begin{array}{*{20}{c}}
I&0\\
0&I
\end{array}} \right]$.

Also, we define $\bar y=\left[ {\begin{array}{*{20}{c}}
y\\
{  {y^*}}
\end{array}} \right]$, $\bar C(a,a^\dagger)=\left[ {\begin{array}{*{20}{c}}
{C(a)}\\
{ C{{(a)}^*}}
\end{array}} \right]$, and $\bar D(t)= \left[ {\begin{array}{*{20}{c}}
{D(t)}&0\\
0&{ D{{(t)}^*}}
\end{array}} \right]$.

Then , the OPO system can be described by
\begin{eqnarray}\nonumber
d\bar a(t)&=& \bar A(a,a^\dagger)(t)dt+\bar B(t)d\bar w(t);\\\nonumber
d\bar y(t)&=& \bar C(a,a^\dagger)(t)dt+\bar D(t)d\bar w(t).
\end{eqnarray}
Note that this OPO system fits within the class of nonlinear quantum systems considered in this paper. In fact, the matrices $A(a,a^\dagger)$ and $ C(a)$ satisfy $\left[C(a), a^T\right]=0$, $\left[A(a,a^\dagger), a^T\right]=-\left[a, A(a,a^\dagger)^T\right]$ and
\begin{eqnarray}\nonumber
&&\frac{1}{{{2{\bar n}}}}\left[ {\bar A{{(a,{a^\dag })}^\dag }{{\bar \theta }^{ - 1}}\bar a,\bar a} \right] - \frac{1}{{{2{\bar n}}}}\left[ {{{\bar a}^\dag }{{\bar \theta }^{ - 1}}\bar A(a,{a^\dag }),\bar a} \right]=\\\nonumber && \qquad \bar A(a,{a^\dag }) - \frac{1}{2}\bar B(a,{a^\dag })\bar C(a,a^\dagger)
\end{eqnarray}
where $\bar n= 3$, $\bar \Theta= \left[ {\begin{array}{*{20}{c}}
I &0\\
0&-I
\end{array}}\right]$ and $\Theta = \left[a,a^\dagger\right]=I$.

In fact, $\frac{1}{{{2{\bar n}}}}\left[ {\bar A{{(a,{a^\dag })}^\dag }{{\bar \theta }^{ - 1}}\bar a,\bar a} \right]= \left[ {\begin{array}{*{20}{c}}
{-\chi ^{(2)}a_1^*a_2}\\
{ \frac{1}{2}\chi ^{(2)} a_1^2}\\
{-\chi ^{(2)}}a_2^{*}{a_1}\\
{ \frac{1}{2} \chi ^{(2)}a_1^{2*}}
\end{array}} \right]$ and $\frac{1}{{{2{\bar n}}}}\left[ {{{\bar a}^\dag }{{\bar \theta }^{ - 1}}\bar A(a,{a^\dag }),\bar a} \right]=\left[ {\begin{array}{*{20}{c}}
{\chi ^{(2)}a_1^*a_2}\\
{ -\frac{1}{2} \chi ^{(2)}a_1^2}\\
{\chi ^{(2)}}a_2^{*}{a_1}\\
{ -\frac{1}{2} \chi ^{(2)}a_1^{2*}}
\end{array}} \right]$. On the other hand, $\bar A(a,a^\dagger)- \frac{1}{2}\bar B(a,{a^\dag })\bar C(a,a^\dagger)=$

$\left[ {\begin{array}{*{20}{c}}
{-2 \chi ^{(2)}}a_1^*a_2\\
{\chi ^{(2)} a_1^2}\\
{-2\chi ^{(2)} a_2^{*}{a_1}}\\
{ \chi ^{(2)} a_1^{2*}}
\end{array}} \right]$. Hence, \begin{eqnarray}\nonumber
&&\frac{1}{{{2{\bar n}}}}\left[ {\bar A{{(a,{a^\dag })}^\dag }{{\bar \theta }^{ - 1}}\bar a,\bar a} \right] - \frac{1}{{{2{\bar n}}}}\left[ {{{\bar a}^\dag }{{\bar \theta }^{ - 1}}\bar A(a,{a^\dag }),\bar a} \right]=\\\nonumber&& \qquad \bar A(a,{a^\dag }) - \frac{1}{2}\bar B(a,{a^\dag })\bar C(a,a^\dagger).
\end{eqnarray}
In addition, the system under consideration is physically realizable. In fact, $\bar B( a, a^\dagger)= \left[\bar C(a,a^\dagger)^\dagger,\bar a \right]\bar J= \left[ {\begin{array}{*{20}{c}}
{ - \sqrt {2{k_{{a_1}}}} }&0&0&0\\
0&{ - \sqrt {2{k_{{a_2}}}} }&0&0\\
0&0&{-\sqrt {2{k_{{a_1}}}} }&0\\
0&0&0&{-\sqrt {2{k_{{a_2}}}} }
\end{array}} \right]$.

On the other hand, $\left[ \bar B( a, a^\dagger),\bar a^\dagger \right] = 0$;
$\left[ \bar a, \bar B( a, a^\dagger)^\dagger \right]=0$ and
$\bar D( a,  a^\dagger)=I$.

Also, $\left[ {\bar A(a,{a^\dag }),{{\bar a}^\dag }} \right] = \left[ {\begin{array}{*{20}{c}}
{ - {k_{{a_1}}}}&{-2 \chi ^{(2)}a_1^*}&{ 2\chi ^{(2)}a_2}&0\\
{ 2 \chi ^{(2)} {a_1}}&{ - {k_{{a_2}}}}&0&{0}\\
{-2\chi ^{(2)}a_2^{*}}&0&{{k_{{a_1}}}}&{ 2 \chi ^{(2)}{a_1}}\\
0&{ 0}&{-2\chi ^{(2)}a_1^*}&{{k_{{a_2}}}}
\end{array}} \right]$,

$\left[ {\bar a,\bar A{{(a,{a^\dag })}^\dag }} \right] =$

$ \left[ {\begin{array}{*{20}{c}}
{ - {k_{{a_1}}}}&{2 \chi ^{(2)}a_1^*}&{ -2\chi ^{(2)}a_2}&0\\
{ -2 \chi ^{(2)} {a_1}}&{ - {k_{{a_2}}}}&0&{0}\\
{2\chi ^{(2)}a_2^{*}}&0&{{k_{{a_1}}}}&{ -2 \chi ^{(2)}{a_1}}\\
0&{ 0}&{2\chi ^{(2)}a_1^*}&{{k_{{a_2}}}}
\end{array}} \right]$ and

$\bar B(a,{a^\dag })\bar J\bar B{(a,{a^\dag })^\dag } = \left[ {\begin{array}{*{20}{c}}
{2{k_{{a_1}}}}&0&0&0\\
0&{2{k_{{a_2}}}}&0&0\\
0&0&{ - 2{k_{{a_1}}}}&0\\
0&0&0&{ - 2{k_{{a_2}}}}
\end{array}} \right]$.

Hence, $\left[ \bar A( a,  a^\dagger),\bar a^\dagger \right]+\left[\bar a, \bar A( a, a^\dagger)^\dagger \right]+ \bar B( a,  a^\dagger)\bar J \bar B( a, a^\dagger)^\dagger =0.$

In this case, the corresponding nonlinear self-adjoint Hamiltonian $\bar H$ is given by
\begin{eqnarray}\nonumber
\bar H&=&\frac{i}{2\bar n}\left({\bar a^\dag }{\bar \theta ^{ - 1}}\bar A(a,{a^\dag }) - \bar A{(a,{a^\dag })^\dag }{\bar \theta ^{ - 1}}\bar a\right)\\\nonumber &=& i\chi ^{(2)}a_1^{2*}a_2-i\chi ^{(2)}a_2^{*}a_1^2.
\end{eqnarray}
and the corresponding coupling operator is $\bar L=\bar C(a,a^\dagger)$.

\section{Conclusion}
In this paper, the issue of physical realizability for a class of nonlinear quantum systems is developed and verified by applying the theory developed to an optical squeezer. A future extension of this paper would be to relate the physical realizability properties to the lossless properties of the class of nonlinear quantum systems considered in this paper.

\end{document}